\documentclass{jnmp} 
\usepackage{amsmath} 
 
\newtheorem*{lemma}{Lemma} 
\newtheorem{proposition}{Proposition} 
 
\theoremstyle{definition} 
\newtheorem*{definition}{Definition} 
\newtheorem{remark}{Remark} 
 
\begin{document} 
 
\renewcommand{\evenhead}{J E Solomin and M Zuccalli} 
\renewcommand{\oddhead}{Noncentral Extensions as Anomalies in 
Classical Dynamical Systems} 
 
\thispagestyle{empty} 
 
\FirstPageHead{10}{1}{2003}{\pageref{solomin-firstpage}--\pageref{solomin-lastpage}}{Letter} 
 
\copyrightnote{2003}{J E Solomin and M Zuccalli} 
 
\Name{Noncentral Extensions as Anomalies\\ in Classical Dynamical 
Systems} 
\label{solomin-firstpage} 
 
\Author{Jorge E SOLOMIN~$^\dag$ and Marcela ZUCCALLI~$^\ddag$} 
 
\Address{$^\dag$~Departamento de Matem\'atica, Facultad de 
Ciencias Exactas,\\ 
~~Universidad Nacional de La Plata and CONICET, Argentina\\ 
~~E-mail: solo@mate.unlp.edu.ar \\[10pt] 
$^\ddag$~Departamento de Matem\'atica, Facultad de Ciencias Exactas,\\ 
~~Universidad Nacional de La Plata, Argentina\\ 
~~E-mail: marce@mate.unlp.edu.ar} 
 
\Date{Received May 14, 2002; Revised 
July 03, 2002; Accepted July 24, 2002} 
 
\begin{abstract} 
\noindent A two cocycle is associated to any action of a Lie group 
on a symplectic manifold. This allows to enlarge the concept of 
anomaly in classical dynamical systems considered by F~Toppan in
[{\it J. 
 Nonlinear Math. Phys.} {\bf 8}, Nr.~3 (2001), 518--533]
so as to encompass some extensions of Lie algebras 
related to noncanonical actions. 
\end{abstract} 
 
\section{Introduction} 
The concept of anomalies in classical dynamical systems was 
considered by F~Toppan in a recent paper \cite{t}. This concept 
allowed him to establish interesting re-interpretations of some 
relevant results. 
 
The general situation he discusses is the following: 
the momentum mappings associated to a canonical action of a 
Lie group $G$  on a classical system  (the Noether 
charges) yield a~nontrivial central extension of ${\rm Lie} (G)$, the 
Lie algebra of $G$. 
 
The aim of this paper is to generalize this kind of 
ideas in order to encompass some noncanonical actions of Lie 
Groups on classical systems that give rise to representations of 
extensions of ${\rm Lie} (G)$. These extensions are in general noncentral. 
 
Our approach is based on the introduction and analysis 
of a two cocycle on ${\rm Lie} (G)$ associated to any action of $G$ on a 
symplectic manifold $(M, \omega)$. This cocycle takes values in 
$C^{\infty}(M)$, but, if the action of $G$ is symplectic, its 
evaluation at any point in $M$ yields a~real valued two cocycle 
which is cohomologous to the cocycle defining the central 
extension mentioned above (see, for instance, \cite{am} or~\cite{w}).

The generalization we propose will allow us to consider 
in our context an example involving the Mickelsson--Fadeev 
cocycle by means of classical objects. This cocycle is associated 
to the quantum anomalous commutator of the constraints of the 
Gauss-law in a $(3+1)$-dimensional Yang--Mills theory interacting 
with Weyl fermions and it was first computed by L~Fadeev by using 
path-integral techniques~\cite{fs} (and also by cohomologi\-cal 
ones~\cite{f}). Almost simultaneously, J~Mickelsson constructed a 
representation of the extension of the current algebra associated 
to this cocycle by using topological techniques. 
 
In fact, we shall see that the Mickelsson--Fadeev 
cocycle is cohomologous, modulo constants, to the cocycle 
associated to the action of the gauge group on the classical 
system defined by the effective lagrangian. 
 
\section{The canonical two-cocycle and the momentum maps} 
 
In this section we consider a symplectic manifold $(M,\omega)$ and an action of a Lie 
group $G$ on $M$ which is not assumed to be symplectic. This action induces a natural 
action of $G$ on the space $C^\infty(M)$: 
 \[ (g\cdot f)(x) = f\left( g^{-1} 
\cdot x\right),\ \ \ \forall g \in G. \] 
 
 The derivative of this action produces a nontrivial action of ${\rm 
Lie}(G)$ on $C^{\infty}(M)$ by \[ a\cdot f = -L_{{\tilde{X}_a}}f, \]with 
${\tilde{X}_a}$ the infinitesimal generator associated to $a\in {\rm Lie} (G)$ and $L$ 
the Lie derivative. 
 
Under such action $C^{\infty}(M)$ becomes a ${\rm Lie}(G)$-module. In general, ${\rm 
Lie}(G)$ acts on all differential forms on $M$ in the same way: $a \cdot \alpha = 
-L_{{\tilde{X}_a}} \alpha$. 
 
 In order to define the cohomology of ${\rm Lie}(G)$ with 
coefficients in $C^{\infty}(M)$, the standard coboundary operator is introduced: 
\begin{gather*} (\delta \alpha)(a_1,\ldots,a_{n+1}) = \sum_{i=0}^{n+1} (-1)^{(i+1)} 
a_i \cdot \alpha (a_1,\ldots,\hat{a_i},\ldots,a_{n+1}) \\ \phantom{(\delta 
\alpha)(a_1,\ldots,a_{n+1}) =}{}+ \sum_{i < j} \alpha ([a_i,a_j],a_1, \ldots, 
\hat{a_i},\ldots,\hat{a_j},\ldots,a_{n+1}) \end{gather*} where $\alpha$ is one $n$ 
cochain on ${\rm Lie}(G)$ with values in $C^{\infty}(M)$ (i.e.\ an alternate 
multilineal map), $a_1,a_2,\ldots,a_{n+1} \in {\rm Lie}(G)$, the symbol `\,\^\,' 
meaning that the variable under it has been deleted and the symbol `\,$\cdot$\,' 
denoting the action of ${\rm Lie}(G)$ on $C^{\infty}(M)$. 
 
The space $Z^n({\rm Lie}(G), C^{\infty}(M))$ of $n$-cocycles consists of 
$n$ cochain $\alpha$ with $\delta \alpha = 0$ and the space 
$B^n({\rm Lie}(G), C^{\infty}(M))$ of $n$-coboundaries consists of $n$ 
cochain such that exists some $(n-1)$ cochain $\beta$ with $\alpha 
= \delta \beta$. 
 
The cohomology groups are defined as $H^n({\rm Lie}(G),C^{\infty}) = 
\displaystyle{\frac {Z^n({\rm Lie}(G), C^{\infty}(M)} {B^n({\rm Lie}(G), 
C^{\infty}(M))}}$. 
 
It is well known that the second cohomology group 
$H^2({\rm Lie}(G),C^{\infty}(M))$ is related to the extensions of 
${\rm Lie}(G)$ by $C^{\infty}(M)$. 
 
The semidirect sum of ${\rm Lie}(G)$ and $C^{\infty}(M)$ consists of pairs $(a,f) \in 
{\rm Lie}(G) \times C^{\infty}(M)$ with the Lie commutator \[ [(a,f),(b,g)] = 
([a,b],\; a\cdot f - b\cdot g). \] 
 
Let $\alpha \in H^2({\rm Lie}(G),C^{\infty}(M))$. We can try to define a modified 
commutator by \[ [(a,f),(b,g)]_{\alpha} = ([a,b],\; a \cdot f - b \cdot g + \alpha 
(a,b)]. \] 
 
(The Jacobi identity for the modifed commutator is easily seen to be equivalent 
to the cocycle condition $\delta \alpha =0$.) 
 
So, each $\alpha \in H^2({\rm Lie}(G),C^{\infty}(M))$ defines a new Lie algebra. 
 
Let $\alpha_1$ and $\alpha_2 \in Z^2({\rm Lie}(G),C^{\infty}( M)).$ 
The Lie algebras formed from these cocycles are isomorphic 
through a mapping of the type $\phi (a,f) = (a, f + \beta(a))$, where 
$\beta \in H^1({\rm Lie}(G), C^{\infty}(M))$. 
 
The condition 
\[ 
[\phi(a,f),\phi(b,g)]_{\alpha_1} = \phi ([(a,f),(b,g)]_{\alpha_2}) 
\] 
is the same as $\alpha_1 - \alpha_2 = \delta \beta$ 
with $\beta \in H^1({\rm Lie}(G),C^{\infty}(M))$ (i.e. $\alpha_1 $ and $\alpha_2$ 
are cohomologous: $\alpha_1 \simeq \alpha_2$). 
 
Thus, up to an isomorphism of the above type 
the Lie algebra extensions are pa\-ra\-met\-ri\-zed 
by elements of $H^2({\rm Lie}(G),C^{\infty}(M))$. 
 
Real valued cocycles and coboundaries can be recovered from the 
previous construction just by considerig real valued cochains as 
constant functions. 
 
Now, we introduce the $C^\infty (M)$-valued two-cocycle on ${\rm Lie} (G)$ canonically 
associated to the action of $G$ mentioned above. 
 
\begin{proposition} Let $\Omega(a,b)(x) = \omega({\tilde{X}_a},{\tilde{X}_b})(x)$ 
$\forall \; a, b\in {\rm Lie}(G)$ and $x \in M$. Then, $\Omega$ is a~$C^\infty 
(M)$-valued two-cocycle on ${\rm Lie}(G)$. \end{proposition} 
 
\begin{proof} 
\begin{gather*} 
(\delta \Omega)(a,b,c) = a \cdot \Omega(b,c) - b \cdot 
\Omega(a,c) + c \cdot \Omega(a,b)\\ 
\phantom{(\delta \Omega)(a,b,c) =}{} - \Omega([a,b],c) + \Omega ([a,c],b) - \Omega ([b,c],a)\\ 
\phantom{(\delta \Omega)(a,b,c)}{}=  - L_{a} \omega ({\tilde{X}_b},{\tilde{X}_c}) 
 + L_{{\tilde{X}_b}} \omega ({\tilde{X}_a},{\tilde{X}_c}) - L_{{\tilde{X}_c}} \omega ({\tilde{X}_a},{\tilde{X}_b})\\ 
\phantom{(\delta \Omega)(a,b,c) =}{}- \omega ({\tilde{X}_c}, \tilde{X}_{[a,b]}) + \omega ({\tilde{X}_b}, \tilde{X}_{[a,c]}) - 
  \omega ({\tilde{X}_a}, \tilde{X}_{[b,c]})\\ 
\phantom{(\delta \Omega)(a,b,c)}{} =  - {\tilde{X}_a} \omega({\tilde{X}_b},{\tilde{X}_c}) + \omega ([{\tilde{X}_a},{\tilde{X}_b}],{\tilde{X}_c}) + \omega ({\tilde{X}_b},[{\tilde{X}_a},{\tilde{X}_c}])\\ 
\phantom{(\delta \Omega)(a,b,c) =}{}+    {\tilde{X}_b} \omega({\tilde{X}_a},{\tilde{X}_c})  - \omega ([{\tilde{X}_b},{\tilde{X}_a}],{\tilde{X}_c}) - \omega ({\tilde{X}_a},[{\tilde{X}_b},{\tilde{X}_c}])\\ 
\phantom{(\delta \Omega)(a,b,c) =}{}-     {\tilde{X}_c} \omega({\tilde{X}_a},{\tilde{X}_b}) + \omega ([{\tilde{X}_c},{\tilde{X}_a}],{\tilde{X}_b})  + \omega ({\tilde{X}_a},[{\tilde{X}_c},{\tilde{X}_b}])\\ 
\phantom{(\delta \Omega)(a,b,c) =}{} -  \omega([{\tilde{X}_a},{\tilde{X}_b}],{\tilde{X}_c}) + \omega([{\tilde{X}_a},{\tilde{X}_c}],{\tilde{X}_b}) - \omega([{\tilde{X}_b},{\tilde{X}_c}],{\tilde{X}_a})\\ 
\phantom{(\delta \Omega)(a,b,c)}{}=   - {\tilde{X}_a} \omega({\tilde{X}_b},{\tilde{X}_c}) +  {\tilde{X}_b} \omega({\tilde{X}_a},{\tilde{X}_c}) -  {\tilde{X}_c} \omega({\tilde{X}_a},{\tilde{X}_b})\\ 
\phantom{(\delta \Omega)(a,b,c) =}{}+ \omega ([{\tilde{X}_a},{\tilde{X}_b}],{\tilde{X}_c}) - \omega ([{\tilde{X}_a},{\tilde{X}_c}],{\tilde{X}_b})  +  \omega ([{\tilde{X}_b},{\tilde{X}_c}],{\tilde{X}_a})\\ 
\phantom{(\delta \Omega)(a,b,c)}{}  =  -3 (d \omega) ({\tilde{X}_a},{\tilde{X}_b},{\tilde{X}_c}) = 0, 
  \end{gather*} 
  where $d$ denotes the usual exterior differential operator. 
\end{proof}

\begin{definition} The cocycle 
\[ 
\Omega(a,b)(x) = \omega ({\tilde{X}_a},{\tilde{X}_b})(x) 
\] 
will be called the {\bf canonical cocycle} associated to the 
action of $G$ on $(M,\omega)$. 
\end{definition} 
 
Thus, each action of $G$ on $M$ yields an extension of ${\rm Lie} (G)$ by 
$C^\infty (M)$. 
 
When the action of $G$ is symplectic, $\Omega (a,b)(x)$ and 
 $\Omega (a,b)(x_0)$ are cohomologous 
as real valued two cocycles on ${\rm Lie} (G)$ for every $x, x_0 \in M$. 
Then, in this case, the element $ [\Omega (a,b)(x)]$ in $H^2({\rm Lie} 
(G), {\mathbb R})$ is independent of $x$. We shall denote it by 
$\tilde{\Omega} (a,b)$. In this way, a central extension of 
${\rm Lie}(G)$ can be defined. 
 
On the other hand, assuming, as henceforth we shall do 
for the symplectic actions to be considered, that the action of 
$G$ admits momentum mappings $J_a$ (that is, for each $ a \in 
{\rm Lie}(G)$ there exists $J_a \in C^{\infty}(M)$ such that $d J_a = 
i_{{\tilde{X}_a}} \omega$),  $\forall \; a,b \in {\rm Lie}(G)$ the function 
 \[ 
\Sigma (a,b)(x):= \{J_a,J_b\}(x) - J_{[a,b]}(x), 
\] 
 with $\{\;,\; \}$  the Poisson bracket 
associated to $\omega$, turns out to be constant on $M$ and, as a~function of $(a,b)$, 
it is a real valued two cocycle on ${\rm Lie} (G)$ \cite{am,w}. This cocycle will be 
denoted by $\tilde {\Sigma} (a,b).$ 
 
Since for any $x_0 \in M$, $J_{[a,b]}(x_0)$ is trivial in 
$H^{2}({\rm Lie} (G), {\mathbb R})$, then 
$\tilde {\Sigma} (a,b)$ is equivalent to $\{J_a,J_b\}(x_0)= \Omega 
(a,b)(x_0)$ in this cohomology. Hence, as real two cocycles on 
${\rm Lie}(G)$, 
\[ 
\{J_a,J_b\}(x) - J_{[a,b]}(x) \simeq \tilde{\Omega} (a,b). 
\] 
 
Thus, for symplectic actions, the momenta give rise to a 
representation of the central extension of ${\rm Lie}(G)$ determined by 
$\tilde{\Omega} (a,b)$. 
 
According to F~Toppan~\cite{t}, an anomaly appears in a 
classical dynamical system when the central extension associated 
to a symplectic action of $G$ on $(M,\omega)$ is nontrivial. Let 
us recall that a necessary condition for it to occur is that 
symplectic potentials are not preserved by the action~\cite{am}. 
 
Now, for a not necessarily symplectic action of $G$, we 
shall show that, under some additional hypotheses, the extension 
of ${\rm Lie}( G)$ defined by $\Omega$ can still be represented by the 
Poisson brackets of some functions associated to the action of~$G$. 
 
More precisely, we shall see that if $\omega$ can be 
written as 
\[ 
\omega = \omega _{i} + \Delta \omega 
\] 
with $\omega_i$ $G$-invariant and $\Delta \omega$ 
a closed form (not necessarily non-degenerate) on $M$ such that 
$L_{{\tilde{X}_a}} \omega = L_{{\tilde{X}_a}} \Delta \omega$ and if $J_a$ are momentum 
maps corresponding to the symplectic action of $G$ on 
$(M,\omega_i)$, we have 
 
\begin{proposition}  Let $\Delta X_a$ denote the vector field 
${\tilde{X}_a} - X_{J_a}$. 
Under the additional assumption 
\[ 
\omega (\Delta X_a, \Delta X_b)=0 \quad  \forall \; a, b \in {\rm Lie}(G) 
\] 
it holds 
\[ 
\{J_a,J_b\} - J_{[a,b]} \simeq \Omega (b,a) \quad \mbox{as} \quad 
C^{\infty}(M)\mbox{-valued cocycle on}\quad  {\rm Lie}(G). 
\] 
\end{proposition} 
 
\begin{proof} 
\begin{gather*} 
  \{J_a,J_b\}   =  \omega(X_{J_a},X_{J_b}) = \omega ({\tilde{X}_a},{\tilde{X}_b}) - 
 \omega({\tilde{X}_a},\Delta X_b) - \omega(\Delta X_a,{\tilde{X}_b})  \\ 
\phantom{\{J_a,J_b\}}{} =  \omega({\tilde{X}_a},{\tilde{X}_b}) + \omega(\Delta X_b,{\tilde{X}_a}) - \omega(\Delta X_a,{\tilde{X}_b}) \\ 
\phantom{\{J_a,J_b\}}{} =  \omega({\tilde{X}_a},{\tilde{X}_b}) + \Delta \omega({\tilde{X}_b},{\tilde{X}_a}) - \Delta \omega({\tilde{X}_a},{\tilde{X}_b})  \\ 
\phantom{\{J_a,J_b\}}{} =  \omega({\tilde{X}_a},{\tilde{X}_b}) - 2\Delta \omega({\tilde{X}_a},{\tilde{X}_b}). 
\end{gather*} 
On the other hand, 
\begin{gather*} 
  (\delta J)(a,b)  =  a\cdot J_b - b\cdot J_a - J_{[a,b]} = 
-L_{\tilde{X}_a} J_b + L_{\tilde{X}_b} J_a - J_{[a,b]} \\ 
\phantom{(\delta J)(a,b)}{} =  -dJ_b({\tilde{X}_a}) + dJ_a({\tilde{X}_b}) - J_{[a,b]}\\ 
\phantom{(\delta J)(a,b)}{} =  -\omega_i({\tilde{X}_b},{\tilde{X}_a}) + \omega_i({\tilde{X}_a},{\tilde{X}_b}) - J_{[a,b]} = 2\omega_i({\tilde{X}_a},{\tilde{X}_b}) - J_{[a,b]}. 
\end{gather*} 
So, $J_{[a,b]} = 2\omega_i({\tilde{X}_a},{\tilde{X}_b}) - (\delta J)(a,b)$. Then, we 
conclude that: 
\begin{gather*} 
  \{J_a,J_b\} -J_{[a,b]}  =  \omega({\tilde{X}_a},{\tilde{X}_b}) - 2\Delta 
\omega({\tilde{X}_a},{\tilde{X}_b}) - 2\omega_i({\tilde{X}_a},{\tilde{X}_b}) + (\delta J)(a,b) \\ 
\phantom{\{J_a,J_b\} -J_{[a,b]} }{}=  \Omega(a,b) +(\delta J)(a,b) 
\end{gather*} 
as we wanted. 
\end{proof} 
 
Thus, under some additional 
hypothesis, noncentral extensions associated to nonsymplectic 
actions can also be represented by physically relevant functions 
and then considered as anomalies of classical dynamical systems. 
 
\begin{remark} The additional hypothesis in the previous proposition is obviously 
fulfilled by any symplectic action since, in this case, $\Delta X_a= 0$ $\forall \; a 
\in {\rm Lie}(G)$. Then, as a $C^\infty (M)$-valued two cocycle on ${\rm Lie}\, (G)$, 
 \[ 
\Sigma (a,b)(x):= \{J_a,J_b\}(x) - J_{[a,b]}(x) \simeq \Omega(b,a)(x). 
\] 
 
On the other hand, as mentioned above, $\Sigma (a,b)(x)$ is 
constant on $M$ and, as real valued cocycle, 
\begin{gather*} 
\tilde{\Sigma} (a,b)= [\Sigma (a,b)(x)] \simeq \Omega(b,a)(x),\\ 
-\tilde{\Sigma} (a,b) = -[\Sigma (a,b)(x)]. 
\end{gather*} 
 
Nevertheless, the opposite signs do not yield a contradiction. In 
fact, from the proof of Proposition~2 we have 
\begin{gather*} 
\tilde{\Sigma} (a,b)= \{J_a,J_b\}(x_0) - J_{[a,b]}(x_0)= 
\Omega(a,b)(x_0) + (\delta J)(a,b)(x_0) \\ 
\phantom{\tilde{\Sigma} (a,b)}{}=\Omega(a,b)(x_0) + 2\omega_i({\tilde{X}_a},{\tilde{X}_b})(x_0) - J_{[a,b]}(x_0)= 
\Omega(a,b)(x_0) - J_{[a,b]}(x_0). 
\end{gather*} 
 
To wit, 
\[ 
\tilde {\Sigma} (a,b)- \omega_i({\tilde{X}_a},{\tilde{X}_b})(x_0) =-J_{[a,b]}(x_0). 
\] 
 
So, in this case, it follows from the Proposition~2 that, as 
expected, 
\[ 
\tilde {\Sigma} (a,b) \simeq \tilde {\Omega}(a,b) 
\] 
as real valued cocycles. 
\end{remark} 
 
\begin{remark} It is shown in \cite{ci} that, under the same hypothesis as in the 
previous proposition, for any $C^{\infty}$-valued two cochain $h(a,b)$ on ${\rm Lie} 
(G)$ such that $dh(a,b) = (\delta \alpha) (a,b)$, $\alpha (a)$ being a local 
symplectic potential of $-L_{\tilde X_a}\omega$, the following equality holds: 
 \[ 
d \{J_a,J_b\} - d J_{[a,b]}  = d h(a,b), 
\] 
and then, 
 \[ 
\{J_a,J_b\} - J_{[a,b]}  =  h(a,b) + c(a,b), 
\] 
 with $c(a,b)$ a real valued two cocycle on ${\rm Lie}(G)$. 
 
 Proposition 2 tells us that, by taking $h(a,b) = \Omega(a,b)+ 
 (\delta J)(a,b)$, we get rid of $c(a,b)$. 
 It is worth to notice that no explicit expression for $h(a,b)$ nor for 
 $c(a,b)$ is given in~\cite{ci}. 
\end{remark} 
 
{\bf A simple example.} Let $G={\mathbb R}^2$ acting by translations on 
$Q={\mathbb R}^2$ and lift the action to $TQ \approx {\mathbb R}^4$. 
 
Let the lagrangian  $L:T {\mathbb R}^2 \simeq {\mathbb R}^4 \rightarrow {\mathbb R}$ 
be defined as 
\[ 
L\big(q^1,q^2,\dot{q^1},\dot{q^2}\big)=\frac {1}{2} \big(\big(\dot{q^1}\big)^2+\big(\dot{q^2} 
\big)^2\big) +{\big(q^1\big)}^2\dot{q^2}. 
\] 
 
It is clear that $L$ can be written as $L=L^i + \Delta L$ with $L^i=\frac {1}{2} 
\big(\big(\dot{q^1}\big)^2+\big(\dot{q^2}\big)^2\big)$ invariant under the action of 
${\mathbb R}^2$ and $\Delta L= {\big(q^1\big)}^2\dot{q^2}$. 
 
The Legendre transform is given by 
\[ 
FL\big(q^1,q^2,\dot{q^1},\dot{q^2}\big)=\big(q^1,q^2,\dot{q^1},\dot{q^1} + 
{(q^1)}^2\big). 
\] 
 
The lagrangian two form  $\omega_L= FL^*\big( dq^1 \wedge dp_1 + dq^2 
\wedge dp_2\big)$ turns out to be 
\begin{gather*} 
\omega_L=\frac {{\partial}^2L}{\partial \dot{q^1} \partial q^2}\, 
d{q^1} \wedge d{q^2} + \frac {\partial^2L}{\partial \dot{q^1} 
\partial \dot{q^2}} \, d{q^1} \wedge d{\dot{q^2}}+ \frac {{\partial}^2L}{\partial \dot{q^2} \partial q^1}\, 
d{\dot{q^1}} \wedge d{q^1} \\ 
\phantom{\omega_L=}{}+ \frac {{\partial}^2L}{\partial 
\dot{q^2} \partial \dot{q^1}} \, d{q^2} \wedge d{\dot{q^1}} 
+ \frac {{\partial}^2L}{\partial \dot{q^1} \partial \dot{q^1}}\, 
d{{q^1}} \wedge d{\dot{q^1}} + \frac 
{{\partial}^2L}{{\partial}\dot{q^2} \partial \dot{q^2}}\, d{q^2} 
\wedge d{\dot{q^2}}. 
\end{gather*} 
 
So, 
\[ 
\omega_L = d{q^1} \wedge d{\dot{q^1}} + d{q^2} \wedge 
d{\dot{q^2}} - 2q^1 d{q^1} \wedge d{q^2}. 
\] 
 
For $a=(a_1,a_2) \in {\mathbb R}^2$, the infinitesimal generator is 
\[ 
{\tilde{X}_a}=(a_1,a_2,0,0)=a_1 \frac {{\partial}}{{\partial}{q^1}} + 
a_2 \frac {{\partial}}{{\partial}{q^2}}. 
\] 
 
The momentum mapping $J_a:{\mathbb R}^2 \rightarrow {\mathbb R}$ must satisfy 
\[ 
dJ_a=i_{{\tilde{X}_a}}\omega_L^i = i_{{\tilde{X}_a}}\omega_L^i=a_1 d\dot{q^1} + a_2 
d\dot{q^2}\] 
and 
\[ 
dJ_a= \frac {\partial J_a}{\partial q^1}\, 
dq^1 + \frac {\partial J_a}{\partial q^2} \, dq^2 + \frac 
{{\partial}J_a}{{\partial}\dot{q^1}} \, d \dot{q^1} + \frac 
{{\partial}J_a}{{\partial}\dot{q^2}} \, d \dot{q^2}. 
\] 
 
We can take 
\[ 
J_a=a_1 \dot{q^1} + a_2 \dot{q^2}. 
\] 
 
The hamiltonian vector associated to the function $J_a$ by 
$\omega_L$ is 
\[ 
X_{J_a}={\tilde{X}_a} - 2q^1 a_2 \frac {{\partial}}{{\partial}\dot{q^1}} - 2q^1 a_1 
\frac {{\partial}}{{\partial}\dot{q^2}}. 
\] 
 
Then, 
\[ 
\Delta X_a=-2q^1\left(a_2 \frac {{\partial}}{{\partial}\dot{q^1}} + 
a_1 \frac {{\partial}}{{\partial}\dot{q^2}}\right). 
\] 
 
Now, it is easy to see that, for this example, it holds \[ \omega_L(\Delta X_a,\Delta 
X_b)=0 \quad \forall \; a, b \in {\rm Lie}\, ({\mathbb R}^2). \] 
 
For $a=(a_1,a_2)$ and $b=(b_1,b_2) \in {\mathbb R}^2$, the 
canonical cocycle turns out to be 
\[ 
\Omega(b,a)=\omega_L({\tilde{X}_b},{\tilde{X}_a})=-2q_1a_1b_2 + 2q^1b_1a_2. 
\] 
 
Thus, 
\[ 
\Omega((1,0),(0,1))=-2q^1. 
\] 
 
On the other hand, a direct computation yields 
\[ 
(\delta J)(a,b) = -L_{{\tilde{X}_a}} {\tilde{X}_b} +L_{{\tilde{X}_b}} {\tilde{X}_a} = -dJ_a ({\tilde{X}_b}) + dJ_b({\tilde{X}_a}) = 0. 
\] 
So, 
\[ 
\{J_a,J_b\}_{\omega _L} = \omega_L (X_{J_a},X_{J_b}) = 
-2q_1a_1b_2 + 2q^1b_1a_2. 
\] 
 
\section{An extension involving the Mickelsson--Faddev cocycle} 
 
The Mickelsson--Fadeev extension of ${\rm Map}\big(S^3,SU(3)\big)$ is a non 
 central one of interest in Quantum Field Theory. It is related to the 
 anomalous commutator of the constraints of the Gauss-law appearing 
 in a $(3+1)$-dimensional Yang--Mills theory interacting 
with Weyl fermions. 
 
 Fifteen years ago, Faddeev~\cite{f} computed this 
anomalous commutator by means of functional integration techniques 
in the following way: 
 
The lagrangian of the theory is given by 
\[ 
{\mathcal L} = tr F^2 + i \bar \psi {\cal D}_A \psi, 
\] 
where $\psi $ is a Weyl fermion, $A \in {\mathcal A}$, the space of 
connections on a trivial bundle over the manifold $S^4$ with 
structure group $SU(3)$, $F=dA$ is it associated curvture 
and ${\mathcal D}_A $ is the covariant 
operator Dirac coupled to~$A$. 
 
The quantization of the theory is carried out in stages. First, the 
fermionic path integral is computed and then, once the effective 
lagrangian is obtained, the bosonic integration is performed. 
 
The first step yields a result which is not gauge invariant: 
\[ 
W[A]= \int D \psi D \bar {\psi}\exp(i \bar \psi {\mathcal D}_A \psi)= \det {\mathcal D}_A . 
\] 
In fact, since ${\mathcal D}_A $ is an unbounded operator, its 
determinant is divergent, a regularizing method  must be applied, 
and this procedure gives rise to a not gauge invariant result. 
 
Let us notice that the variation of $W[A]$ coincides with the variation of the 
Wess--Zumino--Witten lagrangian and is the integrated anomaly of the theory: \[ W[A^g] 
= \exp(i {\alpha}_{1}(A;g)) W[A], \] where $g \in {\mathcal G}$, the gauge group, and 
${\alpha}_{1}(A;g)$ is a 1-cocycle on ${\mathcal G}$ with values in $C^{\infty} 
(\mathcal A)$. 
 
Now, for the second step, the lagrangian that must be considered 
is the effective one: 
\[ 
{\mathcal L}_{E} (A) = \frac {1}{2}\,{\rm tr}\, F_{\mu \nu}F^{\mu \nu} + W[A]. 
\] 
 
Now, the constraints of the Gauss-law are just the momenta (i.e.\ 
the Noether charges) associated to the action of the gauge group 
with respect to the original lagrangian. 
 
By using the Johnson--Bjrken--Low method~\cite{fsb}, Fadeev obtained 
the noncentral extension of ${\rm Lie}(G)$ represented by the 
equal-time commutators of the operators $\tilde {\boldsymbol{G}}_a$, the 
quantum representants of the 
 momenta ${\boldsymbol{G}}_a$: 
\[ [\tilde {\boldsymbol{G}}_a(x),\tilde {\boldsymbol{G}}_b(y)] [A] = 
 C^c_{ab} \tilde {\boldsymbol{G}}_c(y) \delta (x-y) 
+  S_{(a,b)}[A](x,y), 
 \] 
where $ S_{(a,b)}[A](x,y) = \delta 
(x-y) \cdot {\rm MF}_{(a,b)}[A]$ .

The two cocycle ${\rm MF}_{(a,b)}[A] = \displaystyle{\frac {1}{24 \pi} \int_{S^3} 
d^{3}x \,{\rm tr}\,(A[da,db])}$ is {\bf the Mickelsson--Faddeev cocycle}. 
 
In order to relate the Mickelsson--Fadeev cocycle to a 
canonical one, we consider the symplectic manifold $(T 
{\mathcal A},\omega_{{\mathcal L}_E})$, where $\omega_{{\mathcal L}_E}$ is the 
lagrangian form associated to the effective lagrangian ${\mathcal L}_E$. 
 
 Notice that $\omega _{{\mathcal L}_E}$ can be written as $\omega _{{\mathcal L}_E} = \omega 
_{\rm can} + {\Delta }{{\omega}_E}$, where $\omega _{\rm can}$ is the 
symplectic and gauge-invariant structure associated to ${{\mathcal 
L}_{\rm YM}}= \frac {1}{2} \,{\rm tr}\, F_{\mu \nu}F^{\mu \nu} $. 
 
Let us consider the momentum map associated to the lift of the action of the gauge 
group ${\mathcal G}$ on ${\cal A}$ \[ g \cdot A = A^g = g^{-1}Ag + g^{-1}dg \] to 
$(T{\mathcal A},\omega _{\rm can})$.

These maps ${\boldsymbol{G}}_a: T{\mathcal A} \rightarrow {\mathbb R}$ are given by 
$d{\boldsymbol{G}}_a = i_{{{\tilde{X}_a}}}\omega _{\rm can}$ $\forall \; a \in {\rm Lie}({\mathcal G})$ 
where ${\tilde{X}_a}$ is the infinitesimal generator associated to $a\in {\rm Lie} 
(\mathcal G)$. 
 
As mentioned above, these maps are the constraints of the 
Gauss-law~\cite{fs}. 
 
Let us consider the cocycle $\Omega _{\omega _{{\mathcal L}_E}}$ defined on 
${\rm Lie} ({\mathcal G})$ with values in the ${\rm Lie}({\mathcal G})$-module 
$C^{\infty}(T{\mathcal A})$ canonically associated to the action de 
${\mathcal G}$ on ${\mathcal A}$. 
 
The following lemma shows that the additional hypothesis is 
fulfilled in this case.

\begin{lemma} If $\Delta X_a = {\tilde{X}_a} - X_{{\boldsymbol{G}}_a}$, 
where $X_{{\boldsymbol{G}}_a}$ is the Hamiltonian vector field of the 
momentum map ${\boldsymbol{G}}_a$ corresponding to $\omega _{{\mathcal L}_E}$, 
\[ 
\omega_{{\mathcal L}_E} (\Delta X_a,\Delta X_b) = 0 
\quad \forall \; a , b  \in {\rm Lie}({\mathcal G}). 
\] 
\end{lemma} 
 
\begin{proof} 
For any form $\omega$ as in Proposition 2 we have 
\[ 
i_{\Delta X_a} \omega (\cdot)  =  i_{{\tilde{X}_a} - X_{{\boldsymbol{G}}_a }}\omega (\cdot)  = 
 i_{{\tilde{X}_a}}  \omega (\cdot) - i_{X_{{\boldsymbol{G}}_a}} \omega (\cdot)  = 
 i_{{\tilde{X}_a}}  \omega (\cdot) - i_{{\tilde{X}_a}} \omega _{\rm can} (\cdot)  = 
 i_{\tilde{X}_a} \Delta \omega(\cdot). 
\] 
 
Since in the example we are considering , $i_{\Delta X_a} \Delta 
\omega _E= 0$ $\forall \; a \in {\rm Lie}({\mathcal G})$ (formula 4.37 in~\cite{i}), then 
\begin{gather*} 
\omega_{{\mathcal L}_E} (\Delta X_a,\Delta X_b) =  \Delta \omega 
_E({\tilde{X}_a},\Delta X_b)  =  0 \quad \forall \; a, b \in {\rm Lie}({\mathcal G}).\tag*{\qed} 
  \end{gather*} 
  \renewcommand{\qed}{} 
\end{proof} 
 
 Now, it follows from Proposition 2 that 
\[ 
\{J_a,J_b\} - J_{[a,b]} \simeq  \Omega (b,a). 
\] 
 
On the other hand, as shown in~\cite{i}, 
\[ 
\{J_a,J_b\} - J_{[a,b]} \simeq {\rm MF}_{ (a,b)}+ c(a,b). 
\] 
Then, modulo constants, ${\rm MF}_{(a,b)}$ turns out to be cohomologous 
to the canonical cocycle $\Omega (b,a)$. 
 
\label{solomin-lastpage} 
 
\end{document}